\documentclass[12pt]{article}
\usepackage{amsmath, amsfonts}
\usepackage{enumerate}
\usepackage{color}
\numberwithin{equation}{section}
\def \no { \noindent}
\newcommand{\Liminf}{\displaystyle \liminf}

\def\B1{B_{1/2}}

\def\Bl1{{\lambda_1}}
\def\Box{\hfill\rule{2.5mm}{2.5mm}}
\def\C{{\cal {C}}}

\def\H{{\cal H}}

\def\R{{\mathbb {R}}}

\def\bl1{{\bar \lambda_1}}

\def\build#1_#2^#3{\mathrel{
\mathop{\kern 0pt#1}\limits_{#2}^{#3}}}

\def\h1{\mathop{\rm H^1_{\rm loc,\rm u}}}

\def\l2{\mathop{\rm L^2_{\rm loc,\rm u}}}

\newtheorem{cor}{Corollary}[section]

\newtheorem{lem}[cor]{Lemma}
\newtheorem{prop}[cor]{Proposition}
\newtheorem {rem}  {Remark}
\newtheorem{propo}{Proposition}
\newtheorem{thm}[propo]{Theorem}

\begin{document}

\title{\bf The blow-up rate for strongly perturbed semilinear wave equations}

\author{M. A. ${\rm Hamza}$\footnote{The first author is partially supported by the ERC Advanced Grant no.291214, BLOWDISOL during his visit to LAGA, Univ P13 in June 2013. }\,\,  and O. ${\rm Saidi}$}

\maketitle

\begin{abstract}
We consider in this paper a large class of perturbed
semilinear wave equation with subconformal power nonlinearity. In particular, we allow log perturbations of the main source.
We derive a Lyapunov functional in similarity variables and use it  to derive the blow-up rate.
Throughout this work, we use some techniques developed for the unperturbed case studied by Merle and Zaag
\cite{kl} together with ideas introduced by Hamza and Zaag in  \cite{MH} for a class of weather
perturbations.
\end{abstract}

\medskip

{\bf Keywords}: Semilinear wave equation, finite time blow-up,
blow-up rate, perturbations.
\medskip

{\bf MSC 2010 Classification}:
35L05, 35B67, 35B20.

\begin{section}{Introduction}
\end{section}

This paper is devoted to the study of blow-up solutions for the following semilinear
wave equation:

\begin{equation}\label{y}
\left\{
  \begin{array}{ll}
   \partial_{t}^2 u=\Delta u+|u|^{p-1}u+f(u)+g(x,t,\nabla u,\partial_{t} u )\\
  (u(x,0),\partial_{t} u(x,0))=(u_{0}(x),u_{1}(x))\,\,\,\,\,\,\,\,\,\,\,\,\,\,\,\,\,\,\,\,\,\,\,\\
  \end{array}
\right.
\end{equation}
where $u(t):\,\, x\in {\R^N} \rightarrow u(x,t) \in {\R}$, $u_{0}(x)\in H^{1}_{loc,u}$ and $u_{1}(x) \in L^2_{loc,u}$.
The space $L^2_{loc,u}$ is the set of all $v\in L^2_{loc}$ such that
$$\|v\|_{L^2_{loc,u}}\equiv \sup_{d\in{\R^N}}(\int_{|x-d|<1} |v(x)|^2 dx )^{\frac{1}{2}}<+\infty ,$$
and the space $ H^{1}_{loc,u}=\{v \mid v,|\nabla v|\in L^{2}_{loc,u}\}.$
       \\We assume that the functions $f$ and $g$ are $\C^1$, with $f:\R \rightarrow \R$ and
$g: \R^{2N+2} \rightarrow \R$ globally lipschitz, satisfying the following conditions:
$$(H_{f})\,\,\,\,\,\,\,\,\,\,\,\,\,\,\,\,\,\,\,\,\,\,\,|f(x)|\leq M\Big(1+\frac{|x|^p}{(\log (2+x^2))^{a}}\Big),\,\,\,\,\,\,\,\,\,\,\,\,\,\,\,\,\,{\rm for }\,\,{\rm all} \,x\in \R\,\,\,\,\,\,\,\,\,\,{\rm with} \,(M > 0,\,\,\,a> 1),$$
$$(H_{g})\,\,\,\,\,\,\,\,\,|g(x,t,v,z)|\leq  M(1+|v|+|z|), \,\,\,\,\,\,\,\,\,\,\,\,\,{\rm for }\,\,{\rm all} \,x,\,v\in \R^N\,t,z\in \R\,{\rm with }\,(M > 0).$$
Finally, we assume that

$$p>1\,\,{\rm and} \,\,p < p_{c}\equiv 1+\frac{4}{N-1}\,\,\,{\rm if} \,\,N\geq 2.$$

We can mention that Merle and Zaag \cite{fh3}, \cite{kl} and \cite{fh4}, studied the equation $\eqref{y}$ with $(f,g)\equiv(0,0)$ and $1<p\leq p_{c}$, later we give more details for this fact.
Then in \cite{MH} and \cite{MH1}, Hamza and Zaag consider a similar class of perturbed equations, with $(H_{f})$ and $(H_{g})$ replaced by a more restrictive conditions: $|f(u)|\leq M(1+|u|^q)$ and $|g(u)|\leq M(1+|u|)$ for some $q<p,\,\,M > 0.$ In this paper, we go beyond this restriction and allow stronger perturbations.

The Cauchy problem of equation $\eqref{y}$ is wellposed in $H^{1}_{loc,u}\times L^{2}_{loc,u}$. This follows from the finite speed of propagation and the wellposdness in $H^{1}\times L^{2}$, valid whenever
$1< p <1+\frac{4}{N-2}$. The existence of blow-up solutions $u(t)$ of $\eqref{y}$ follows from energy techniques (see
 for example Levine and Todorova \cite{HAG}, and Todorova \cite{GT}).
 If $u(t)$ is a blow-up solution of
$\eqref{y}$, we define (see for example Alinhac \cite{sa} and \cite{sa1}) $\Gamma$ as the graph of a function $x \mapsto T(x)$
 such that the domain of definition of $u$ (also called the maximal influence domain)
$$D_{u}=\{ (x,t)| t< T(x)\}.$$
Moreover, from the finite speed of propagation, $T$ is a 1-Lipschitz function.
Let us first introduce the following non degeneracy condition for $\Gamma$. If we introduce
  for all $x \in \R^N$, $t \leq T(x)$ and $\delta > 0$, the cone

  \begin{equation}
   C_{x,t,\delta}=\{(\xi,\tau)\neq (x,t)| 0\leq \tau \leq t- \delta | \xi-x|\},\\
 \end{equation}
 then our non degeneracy condition is the following: $x_{0} $ is a non-characteristic point if

\begin{equation} \label{17}
 \exists \delta_{0}=\delta_{0}(x_{0})\in (0\,\,1)\,\,{\rm such }\,\,{\rm that} \,\,u\,\, {\rm is}\,\, {\rm defined}\,\, {\rm on}\,\, C_{x_{0},T(x_{0}),\delta_{0}}.\\
 \end{equation}
Hamza and Zaag in \cite{MH} and \cite{MH1} have proved, that if $u$ is a solution of $\eqref{y}$ with $(H_{f})$ and $(H_{g})$ replaced by the following condition: $|f(u)|\leq M(1+|u|^q)$ and $|g(u)|\leq M(1+|u|)$ for some $q<p,\,\,M > 0,$ with blow-up graph $\Gamma : \{ x\mapsto T(x)\}$ and $x_{0}$ is a non-characteristic
  point (in the sense $\eqref{17}$), then there exist $S_{1}>0$ such that, for all $t\in [t_{2}(x_{0}),T(x_{0}))$, where $t_{2}(x_{0})=T(x_{0})-e^{-s_{2}(x_{0})} $ and $s_{2}(x_{0})=\max(S_{1},-\log(\frac{T(x_{0})}{4}))$, we have
 \begin{eqnarray}\label{28nov}
  0 < \varepsilon_{0}(N,p)&\leq &(T(x_{0})-t)^{\frac{2}{ p-1}}\frac{\|u(t)\|_{L^{2}(B(x_{0},T(x_{0})-t))}}{(T(x_{0})-t)^{\frac{N}{2}}}\nonumber\\
  &&+(T(x_{0})-t)^{\frac{2}{ p-1}+1}\Big(\frac{\|\partial_{t} u(t)\|_{L^{2}(B(x_{0},T(x_{0})-t))}}{(T(x_{0})-t)^{\frac{N}{2}}}\\
  &&+
 \frac{\|\nabla u(t)\|_{L^{2}(B(x_{0},T(x_{0})-t))}}{(T(x_{0})-t)^{\frac{N}{2}}}\Big)\leq K,\nonumber
  \end{eqnarray}
 where $K=K(s_{2}(x_{0}),\|(u(t_{2}(x_{0})),\partial_{t} u(t_{2}(x_{0})))\|_{H^{1}\times L^{2}(B(x_{0},\frac{e^{-s_{2}(x_{0})}}{\delta_{0}(x_{0})}))})$ and the constant $\delta_{0}(x_{0})\in (0,1)$ is defined in $\eqref{17}$.
 \\In the pure power case, equation $\eqref{y}$ reduces to the semilinear wave equation:

  \begin{equation} \label{16}
  \partial_{t}^2 u =\Delta u+|u|^{p-1}u,\,\,\,(x,t)\in \R^N \times [0,T).\\
  \end{equation}
\bigskip

It is interesting to recall that previously Merle and Zaag in \cite{fh3} and \cite{kl}  proved a similar result as \eqref{28nov} for  $u$ a solution of $\eqref{16}$ with blow-up graph $\Gamma : \{ x\mapsto T(x)\}.$ 
\\We aim for studying the growth estimate for $u(t)$ near the space time blow-up graph, and extend the result of Hamza and Zaag \cite{MH} for a larger class of perturbation.
\\In order to keep our analysis clear, we may assume that $f(u)\equiv \frac{|u|^{p}}{(\log (2+u^2))^{a}}$ and $g \equiv 0$, in the equation $\eqref{y}$. The adaptation to the case $g\not\equiv  0$ is straightforward  from the techniques.

\bigskip

As in \cite {MH}, \cite {fh3}, \cite {kl} and \cite {fh4}, we want to write the solution $v$ of the associate ordinary differential equation of
   $\eqref{y}$. It is clear that $v$ is given by
   \begin{equation}
   v"=v^p+ f(v),\,\,\,v(T)=+\infty,\\
   \end{equation}
   and satisfies: $v(t)\sim \frac{\kappa}{(T-t)^{\frac{2}{p-1}}}$ as $t\longrightarrow T$, where $\kappa = \Big(\frac{2p+2}{(p-1)^2}\Big)^{\frac{1}{p-1}}.$
 For this reason, we define for all $ x_{0} \in \R^N$, $0< T_{0} \leq T_{0}(x_{0})$, the following similar transformation
  introduced in Antonini and Merle \cite {cf} and used in \cite {MH}, \cite {MH1}, \cite {fh3}, \cite {kl} and \cite {fh4}:

  \begin{equation} \label{2}
y=\frac{x-x_{0}}{T_{0}-t}\,\,\,\,\,\,s=-\log(T_{0}-t),\,\,{\rm and }\,\, w_{x_{0},T_{0}}(y,s)=(T_{0}-t)^\frac{2}{p-1}u(x,t).\,\,\,\,\,\,\,\,\,\,\,\\
\end{equation}
The function $w_{x_{0},T_{0}}$  (we write $w$ for simplicity) satisfies the following equation for all
 $y\in B$ and $s\geq -\log(T_{0})$:

\begin{eqnarray}\label{1}
\partial_{s}^2w&=&\frac{1}{\rho} div(\rho\nabla w-\rho (y.\nabla w)y)-\frac{2(p+1)}{(p-1)^2}w+|w|^{p-1}w\nonumber \\
&&-\frac{p+3}{p-1}\partial_{s}w-2(y.\nabla \partial_{s}w)+e^{\frac{-2ps}{p-1}}f(e^{\frac{2s}{p-1}}w),
\end{eqnarray}
where  $ \rho(y)=(1-|y|^2)^{\alpha}$ and $\alpha=\frac{2}{p-1}- \frac{N-1}{2}>0$.
In the new set of variable $(y,s)$, the behavior of $u$ as $t\rightarrow T_{0}$ is equivalent
to the behavior of $w$ as $s\rightarrow +\infty.$ The equation $\eqref{1}$ will be studied in the space $\H$

 $$\H =\{(w_{1},w_{2})| \int_{B} \Big(w_{2}^2+|\nabla w_{1}|^2(1-|y|^2)+w_{1}^2\Big)\rho dy< +\infty\}. $$
In the whole paper we denote
 \begin{equation}\label{90}
 F(u)=\int_{0}^{u}f(v)dv.\\
 \end{equation}
In the non-perturbed case, Antonini and Merle \cite {cf} proved that
\begin{equation}\label{110}
E_{0}(w(s))=\int_{B}(\frac{1}{2}(\partial_{s}w)^2+\frac{1}{2}|\nabla w|^2-\frac{1}{2}(y.\nabla w)^2+\frac{p+1}{(p-1)^2}w^2-\frac{1}{p+1}|w|^{p+1})\rho dy,\,\,\,\,\,\,\,\,\,\,\,\,\,\,\,\,\\
\end{equation}
is a Lyapunov functional for equation $\eqref{1}$. In our case, we introduce

\begin{equation}\label{3}
H(w(s),s)=\exp\Big(\frac{p+3}{(a-1)s^{b-1}}\Big) E(w(s),s)+\theta e^{\frac{-(p+1)s}{p-1}},{\rm with } \,\,b=\frac{a+1}{2}\,\,\,\,\,\,\,\,\,\\
\end{equation}
where $\theta$ is a sufficiently large constant that will be determined later,

\begin{eqnarray}\label{4}
E(w(s),s)&=&E_{0}(w(s))+I(w(s),s)+J(w(s),s),\\
I(w(s),s)&=&-e^{\frac{-2(p+1)s}{p-1}} \int_{B} F(e^{\frac{2s}{p-1}}w)\rho dy,\nonumber\\
J(w(s),s)&=&-\frac{1}{s^{b}}\int_{B} w\partial_{s}w \rho dy.\nonumber
\end{eqnarray}
We now claim that the functional $H(w(s),s)$ is a decreasing function of time for equation $\eqref{1}$,
provided that $s$ is large enough.
\\Here we announce our main result.

 \begin{thm}

Let N, p, a and M be fixed. There exists $S_{1}=S_{1}(N,p,a,M)\geq 1 $ such that, for all $s_{0}\in \R$ and $w$ solution of equation $\eqref{1}$ satisfying
 $(w,\partial_{s}w)\in C([s_{0},\,\,+\infty),\H)$, it holds that $H(w(s),s)$ satisfies the following inequality, for all $s_{2}> s_{1} \geq \max(S_{1},s_{0})$,

\begin{equation}
 H(w(s_{2}),s_{2})-H(w(s_{1}),s_{1}) \leq -\alpha \int_{s_{1}}^{s_{2}}\int_{B} (\partial_{s}w)^2\frac{\rho}{1-|y|^2}dyds.\\
\end{equation}
\end{thm}

\newpage
 \begin{rem}
\quad
{\rm  \begin{enumerate}[(i)]
\item Our method breaks down in the conformal case when $p\equiv p_{c}$, since in the
energy estimates in similarity variables, the perturbations terms are integrated
on the whole unit ball, hence, difficult to control with the dissipation of the non perturbed
equation  $\eqref{16}$, which degenerates to the boundary of the unit ball.
  \item  The existence of this Lyapunov functional (and a blow-up criterion
  for equation  $\eqref{1}$ based in $H(w(s),s)$, see Lemma 3.1 below) are a crucial step in the derivation
  of the blow-up rate for equation $\eqref{y}$. Indeed, with the functional $H(w(s),s)$ and some more work,
  as in \cite {MH} and \cite {MH1} we are able to adapt the analysis performed in \cite {kl} for equation $\eqref{16}$ and get the Theorem 2 below.
 \item It is worth noticing that the method breaks down when $a\leq 1$ too,
because with some analysis to the Lyapunov functional we find an equality of type $\frac{d}{ds}(E(w(s))\leq \frac{C}{s^a}E(w(s))$, $E$ is upper bounded if $a>1$ but if $a\leq 1$ we can not conclude.

\end{enumerate}}

\end{rem}

  \begin{thm}

Let N, p, a and M be fixed. There exists $\hat{S_{0}}= \hat{S_{0}}(N,p,a,M)\in \R $, and $\varepsilon_{0}=\varepsilon_{0}(N,p,a,M)$, such that if u is a solution of
$\eqref{y}$ with blow-up graph $\Gamma: \{x\mapsto T(x)\}$ and $x_{0}$ is a non-characteristic point, then
\begin{enumerate}[{\rm (i)}]
\item For all  $s \geq \hat{s_{0}}(x_{0})=\max(\hat{S_{0}}(N,p,a,M),-\log(\frac{T(x_{0})}{4}))$,
$$0 < \varepsilon_{0} \leq \|w_{x_{0},T(x_{0})}(s)\|_{H^{1}(B)} + \|\partial_{s}w_{x_{0},T(x_{0})}(s)\|_{L^{2}(B)}\leq K.$$

\item  For all $t\in [t_{0}(x_{0}),T(x_{0}))$, where $t_{0}(x_{0})=\max(T(x_{0})-e^{-\hat{s_{0}}(x_{0})},\frac{3T(x_{0})}{4})$, we have
\begin{eqnarray*}
0 &< &\varepsilon_{0}(N,p)\leq (T(x_{0})-t)^{\frac{2}{ p-1}}\frac{\|u(t)\|_{L^{2}(B(x_{0},T(x_{0})-t))}}{(T(x_{0})-t)^{\frac{N}{2}}}\\
&&+(T(x_{0})-t)^{\frac{2}{ p-1}+1}\Big(\frac{\|\partial_{t}u(t)\|_{L^{2}(B(x_{0},T(x_{0})-t))}}{(T(x_{0})-t)^{\frac{N}{2}}}\\
&&+\frac{\|\nabla u(t)\|_{L^{2}(B(x_{0},T(x_{0})-t))}}{(T(x_{0})-t)^{\frac{N}{2}}}\Big)\leq K,
\end{eqnarray*}
\end{enumerate}
where $K=K(N,p,a,\hat{s_{0}}(x_{0}),\parallel(u(t_{0}(x_{0})),\partial_{t}u(t_{0}(x_{0}))\parallel_{H^{1}\times L^{2}(B(x_{0},\frac{e^{-\hat{s_{0}}(x_{0})}}{\delta_{0}(x_{0})}))}$ and $\delta_{0}(x_{0})\in (0,1)$ is defined in $\eqref{17}$.
\end{thm}

\begin{rem}
{\rm  In a series of papers \cite {fh3}, \cite {kl}, \cite {fh4}, \cite {fh}, \cite {fh1} and \cite {fh2}, Merle and Zaag give a full picture of the blow-up for solution of $\eqref{16}$, in one space dimension and in dimension space $N\geq 2$.
\\The result of all this paper is extended by Hamza and Zaag for a class of perturbed problem in one space dimension or in higher dimension under radial symetry outside origin in \cite {MH2}, or in dimension space $N\geq 2$ in  \cite{MH} and \cite {MH1}, (blow-up, profile, regularity of the blow-up graph, existence
 of characteristic points, etc...). Once again, we believe that the key point in the analysis of blow-up for equation $\eqref{y}$ is the derivation of a Lyapunov
 functional in similarity variables, which is the object of our paper.
 \\\\As in \cite{MH}, \cite {MH1}, \cite {MH2}, \cite {fh3}, \cite {kl}, \cite {fh4}, \cite {fh}, \cite {fh1} and \cite {fh2}, the proof of Theorem 2 relies on four ideas (the existence of a Lyapunov functional, interpolation
 in Sobolev spaces, some Gagliardo-Nirenberg estimates and a covering technique adapted to the geometric shape of the blow-up surface).
  It happens that adapting the proof of \cite {fh3}, \cite {kl} and \cite {fh4} given in the non-perturbed case $\eqref{16}$ is straightforward, except for a key argument, where we
  bound the $L^{p+1}$ space-time norm of $w$. Therefore, we only present that argument, and refer to \cite {fh3}, \cite {kl} and \cite {fh4}, for the rest of the proof.}
 \end{rem}
This paper is divided in two sections, each of them devoted to the proof of a Theorem.

\begin{section}{A Lyapunov functional for equation $\eqref{1}$}
\end{section}

\medskip

Throughout this section, we prove Theorem 1, we consider $(w,\partial_{s} w)\in C([s_{0},\,\,+\infty),\H)$ where $w$ is a solution
of $\eqref{1}$ and $s_{0}\in \R $. We aim at proving that the functional $H(w(s),s)$ defined in $\eqref{3}$ is a Lyapunov functional for equation $\eqref{1}$, provided that $s \geq S_{1}$, for some $S_{1}=S_{1}(N,p,M,a)$.
We denote the unit ball of $\R^N$ by $B$. We denote by $C$ a constant which depends only on $N,p,a$ and on $|B|$.
\\The starting point in our analysis is to prove the following lemma.

\begin{lem}

 For all $s\geq \max(s_{0},1)$,

  \begin{equation}\label{300}
\frac{d}{ds}(E_{0}(w(s))+I(w(s),s))= -2\alpha \int_{B}(\partial_{s}w)^2 \frac{\rho}{1-|y|^2}dy +\Sigma_{0}(s),\,\,\,\,\,\,\,\,\,\,\,\,\,\,\,\,\\
\end{equation}
{\rm where} $\Sigma_{0}(s)$ {\rm satisfies}

 \begin{equation}\label{301}
\Sigma_{0}(s)\leq  Ce^{\frac{-(p+1)s}{p-1}}+ \frac{C}{s^{a}}\int_{B}|w|^{p+1}\rho dy
.\\
\end{equation}

\end{lem}

 \bigskip

{\it Proof}: Multiplying $\eqref{1}$  by $\partial_{s}w \rho$ and integrating over the unit ball $B$, we obtain, for all
   $s\geq s_{0}$,

\begin{equation}\label{200}
\frac{d}{ds}(E_{0}(w(s))+I(w(s),s))=-2\alpha \int_{B}(\partial_{s}w)^2\frac{\rho}{1-|y|^2}dy+\Sigma_{0}^{1}(s)+ \Sigma_{0}^{2}(s),\\
\end{equation}
where
\begin{equation}\label{XX}
\Sigma_{0}^{1}(s)=\frac{2(p+1)}{p-1}e^{\frac{-2(p+1)s}{p-1}}\int_{B} F(e^{\frac{2s}{p-1}}w)\rho dy \,\,{\rm and}\,\, \Sigma_{0}^{2}(s)=-\frac{2e^{\frac{-2ps}{p-1}}}{p-1}\int_{B}f(e^{\frac{2s}{p-1}}w)w\rho dy.\\
\end{equation}
It is clear that we obtain $\eqref{300}$ with $\Sigma_{0}(s)=\Sigma_{0}^{1}(s)+\Sigma_{0}^{2}(s)$. Now, in order  to obtain estimate $\eqref{301}$, it is enough to control the terms $\Sigma_{0}^{1}(s)$ and $\Sigma_{0}^{2}(s)$.
\\Clearly the function $F$ defined in $\eqref{90}$ satisfies the following estimate:
\begin{equation}\label{91}
|F(x)|+|xf(x)|\leq C(1+\frac{ |x|^{p+1}}{(\log(2+x^2))^a}).\\
\end{equation}
Taking advantage of inequality $\eqref{91}$, we see that
\begin{equation}\label{101}
 |\Sigma_{0}^{1}(s)|+|\Sigma_{0}^{2}(s)| \leq Ce^{\frac{-2(p+1)s}{p-1}}+C \int_{B}\frac{|w|^{p+1}}{(\log(2+e^{\frac{4s}{p-1}}w^2))^a}\rho dy.\\
\end{equation}
In order to prove $\eqref{301}$, we divide the unit ball $B$ in two parts
 $$A_{1}(s)=\{y \in B\,\,|\,\, w^2(y,s)\leq  e^{\frac{-2s}{p-1}}\}\,\,{\rm and }\,\,A_{2}(s)=\{y \in B\,\,|\,\, w^2(y,s)>  e^{\frac{-2s}{p-1}}\}.$$
It follows then that
\begin{equation}\label{103}
\int_{B}\frac{|w|^{p+1}}{(\log(2+e^{\frac{4s}{p-1}}w^2))^a}\rho dy=\int_{A_{1}(s)}\frac{|w|^{p+1}}{(\log(2+e^{\frac{4s}{p-1}}w^2))^a}\rho dy+\int_{A_{2}(s)}\frac{|w|^{p+1}}{(\log(2+e^{\frac{4s}{p-1}}w^2))^a}\rho dy.
\end{equation}
On the one hand, ${\rm if }\,\,y \in A_{1}(s)$, we have
$$\frac{|w|^{p+1}}{(\log(2+e^{\frac{4s}{p-1}}w^2))^a} \leq \frac{e^{\frac{-(p+1)s}{p-1}}}{(\log(2))^a}.$$
If we integrate over $A_{1}(s)$,
using the fact that $\int_{A_{1}(s)}\rho dy\leq |B|$, we see that
\begin{equation}\label{93}
\int_{A_{1}(s)}\frac{|w|^{p+1}}{(\log(2+e^{\frac{4s}{p-1}}w^2))^a}\rho dy \leq \frac{e^{\frac{-(p+1)s}{p-1}}}{(\log(2))^a}\int_{A_{1}(s)}\rho dy\leq C e^{\frac{-(p+1)s}{p-1}}. \\
\end{equation}
On the other hand, if $y\in A_{2}(s)$, we have
$$\log(2+e^{\frac{4s}{p-1}}w^2)>\log(2+e^{\frac{2s}{p-1}}) \geq \frac{2s}{p-1},  $$
and for all $s\geq \max (s_{0},1)$, we write for all$\,\,y \in A_{2}(s)$,
$$\frac{|w|^{p+1}}{(\log(2+e^{\frac{4s}{p-1}}w^2))^a} \leq \frac{C}{s^a} |w|^{p+1}.$$
We integrate now over $A_{2}(s)$, using the simple fact that $A_{2}(s)\subset B,$ we obtain for all $s\geq \max (s_{0},1)$,
\begin{equation}\label{102}
\int_{A_{2}(s)}\frac{|w|^{p+1}}{(\log(2+e^{\frac{4s}{p-1}}w^2))^a}\rho dy\leq \frac{C}{s^a}\int_{A_{2}(s)}|w|^{p+1} \rho dy\leq \frac{C}{s^a}\int_{B}|w|^{p+1} \rho dy.\\
\end{equation}
To conclude, it suffices to combine $\eqref{101}$, $\eqref{103}$, $\eqref{93}$ and $\eqref{102}$, then write
\begin{equation}\label{150}
|\Sigma_{0}^{1}(s)|+|\Sigma_{0}^{2}(s)| \leq Ce^{\frac{-(p+1)s}{p-1}} +\frac{C}{s^a}\int_{B}|w|^{p+1} \rho dy,\\
\end{equation}
which ends the proof of Lemma 2.1.

\Box
 \bigskip

  We are going now to prove the following estimate for the functional $J(w(s),s)$:

\begin{lem}

For all $s\geq \max(s_{0},1)$,

\begin{eqnarray}
\frac{d}{ds}(J(w(s),s))&\leq &\alpha\int_{B}(\partial_{s}w)^2\frac{\rho}{1-|y|^2}dy +\frac{p+3}{2s^{b}}E((w(s),s)) \nonumber\\
                            &&-\frac{p-1}{4s^{b}}\int_{B}|\nabla w|^2(1-|y|^2)\rho dy-\frac{p+1}{2(p-1)s^{b}}\int_{B} w^2 \rho dy \nonumber\\
 &&-\frac{p-1}{2(p+1)s^{b}}\int_{B}|w|^{p+1}\rho dy+\Sigma_{1}(s),
 \end{eqnarray}
{\rm where} $\Sigma_{1}(s)$ {\rm satisfies the following inequality:}

\begin{eqnarray}
\Sigma_{1}(s)&\leq & \frac{C}{s^{2b}}\int_{B} w^2 \rho dy+ \frac{C}{s^{2b}}\int_{B} |\nabla w|^2(1-|y|^2)\rho dy\nonumber\\
&&+ \frac{C}{s^{a+b}}\int_{B} |w|^{p+1} \rho dy+ C e^{\frac{-(p+1)s}{p-1}}.
\end{eqnarray}
\end{lem}

{\it Proof}: Note that $J(w(s),s)$ is a differentiable function for all $s\geq s_{0}$ and that
  $$  \frac{d}{ds}(J(w(s),s))=\frac{b}{s^{b+1}}\int_{B} w\partial_{s}w \rho dy-\frac{1}{s^{b}}\int_{B} (\partial_{s}w)^2 \rho dy -\frac{1}{s^{b}}\int_{B} w\partial_{s}^2 w \rho dy.$$
By using the equation $\eqref{1}$ and integrating by parts, we have
\begin{eqnarray}\label{111}
\frac{d}{ds}(J(w(s),s))&=&-\frac{1}{s^{b}}\int_{B} (\partial_{s}w)^2 \rho dy+\frac{1}{s^{b}}\int_{B}(|\nabla w|^2-(y.\nabla w)^2)\rho dy+ \frac{2(p+1)}{(p-1)^2s^{b} }\int_{B} w^2 \rho dy\nonumber \\
&&-\frac{1}{s^{b}}\int_{B}|w|^{p+1}\rho dy+\Sigma_{1}^{1}(s)+\Sigma_{1}^{2}(s)+\Sigma_{1}^{3}(s)+\Sigma_{1}^{4}(s),
 \end{eqnarray}

 where
 \begin{eqnarray*}
\Sigma_{1}^{1}(s)&=&(\frac{b}{s}+\frac{p+3}{p-1}-2N) \frac{1}{s^{b}}\int_{B} w\partial_{s}w \rho dy\\
\Sigma_{1}^{2}(s)&=&-\frac{2}{s^{b}}\int_{B}\partial_{s}w(y.\nabla w)\rho dy\\
\Sigma_{1}^{3}(s)&=&-\frac{e^{\frac{-2ps}{p-1}}}{s^{b}}\int_{B} wf(e^{\frac{2s}{p-1}}w)\rho dy\\
\Sigma_{1}^{4}(s)&=&-\frac{2}{s^{b}}\int_{B}w\partial_{s}w(y.\nabla \rho)dy.
\end{eqnarray*}
By combining $\eqref{110}$, $\eqref{4}$, $\eqref{111}$ and some straightforward computations,  we see that
\begin{eqnarray}\label{112}
\frac{d}{ds}(J(w(s),s))&=& -\frac{p+7}{4s^{b}}\int_{B}(\partial_{s}w)^2 \rho dy +\frac{p+3}{2s^{b}}E(w(s),s) -\frac{p-1}{4s^{b}}\int_{B}(|\nabla w|^2-(y.\nabla w)^2)\rho dy\nonumber\\
 &&-\frac{p+1}{2(p-1)s^{b}}\int_{B}w^2 \rho dy-\frac{p-1}{2(p+1)s^{b}}\int_{B}|w|^{p+1}\rho dy\nonumber\\
&&+\Sigma_{1}^{1}(s)+\Sigma_{1}^{2}(s)+\Sigma_{1}^{3}(s)+\Sigma_{1}^{4}(s)+\Sigma_{1}^{5}(s)+\Sigma_{1}^{6}(s),
\end{eqnarray}
  where
\begin{eqnarray*}
\Sigma_{1}^{5}(s)&=&\frac{p+3}{2s^{b}}\int_{B} w\partial_{s}w \rho dy\\
\Sigma_{1}^{6}(s)&=&\frac{p+3}{2s^{b}}e^{\frac{-2(p+1)s}{p-1}} \int_{B} F(e^{\frac{2s}{p-1}}w)\rho dy.
\end{eqnarray*}
We now study each of the last six terms. To estimate $\Sigma_{1}^{1}(s)$ and $\Sigma_{1}^{5}(s)$, using the fact that for all $s\geq\max(s_{0},1)$,
  $$\Big|\frac{b}{s}+\frac{p+3}{p-1}-2N+\frac{p+3}{2s^{b}}\Big|\leq C,$$
   we get by virtue of Cauchy-Schwartz inequality:
    \begin{equation}\label{6}
   |\Sigma_{1}^{1}(s)|+|\Sigma_{1}^{5}(s)|\leq  \frac{C}{s^{b}}\int_{B}|w\partial_{s}w| \rho dy\leq \frac{\alpha}{3}\int_{B}(\partial_{s}w)^2\frac{\rho}{1-|y|^2}dy+ \frac{C}{s^{2b}}
  \int_{B}w^2\rho dy.\\
  \end{equation}
  Using again the Cauchy-Schwartz inequality, we obtain
  \begin{equation}\label{10}
  |\Sigma_{1}^{2}(s)|\leq  \frac{\alpha }{3}\int_{B}(\partial_{s}w)^2\frac{\rho}{1-|y|^2}dy+\frac{C}{s^{2b}}\int_{B} |\nabla w|^2 \rho (1-|y|^2)dy.\\
  \end{equation}
Using $\eqref{XX}$, we write for all $s\geq \max(s_{0},1)$,
  $$\Sigma_{1}^{3}(s)=\frac{(p-1)}{2s^b}\Sigma_{0}^{2}(s),\,\,{\rm and }\,\,\Sigma_{1}^{6}(s)=\frac{(p+3)(p-1)}{4(p+1)s^b}\Sigma_{0}^{1}(s).$$
This easily leads to the following result
$$|\Sigma_{1}^{3}(s)|+ |\Sigma_{1}^{6}(s)|\leq \frac{C}{s^b}(|\Sigma_{0}^{1}(s)|+|\Sigma_{0}^{2}(s)|).$$
By exploiting inequality $\eqref{150}$ and the fact that $s\geq 1$, we see that
  \begin{equation}\label{12}
 |\Sigma_{1}^{3}(s)|+ |\Sigma_{1}^{6}(s)|\leq C e^{\frac{-(p+1)s}{p-1}}+\frac{C}{s^{a+b}}\int_{B}|w|^{p+1}\rho dy .\,\,\,\,\,\,\,\,\,\,\,\,\,\\
   \end{equation}
Now, we estimate the expression $\Sigma_{1}^{4}(s)$. Since we know that $y.\nabla \rho=-2\alpha\frac{|y|^2}{(1-|y|^2)}\rho$, we can use the
Cauchy-Schwartz inequality to write
 $$ |\Sigma_{1}^{4}(s)|\leq  \frac{C}{s^{b}}\int_{B}|\partial_{s}w|(1-|y|^2)^{\frac{\alpha -1}{2}}|w||y|(1-|y|^2)^{\frac{\alpha -1}{2}}dy,$$
   \begin{equation}\label{7}
  \,\,\,\,\,\,\,\,\,\,\,\,\,\,\,\,\,\,\,\,\leq \frac{\alpha }{3}\int_{B}(\partial_{s}w)^2\frac{\rho}{1-|y|^2}dy+ \frac{C}{s^{2b}} \int_{B}w^2\frac{|y|^2\rho}{1-|y|^2}dy.\,\,\,\\
  \end{equation}
 Since we have the following Hardy type inequality for any  $w \in H^{1}_{loc,u}({\R^N})$ (for more details on this subject, we
refer the reader to appendix B in \cite{fh3}):
 \begin{equation}\label{8}
  \int_{B}w^2\frac{|y|^2\rho}{1-|y|^2}dy\leq C \int_{B} |\nabla w|^2 \rho (1-|y|^2)dy+C\int_{B}w^2 \rho dy,\\
  \end{equation}
we get from $\eqref{7}$ and $\eqref{8}$
    \begin{equation}\label{9}
  |\Sigma_{1}^{4}(s)|\leq  \frac{\alpha }{3}\int_{B}(\partial_{s}w)^2\frac{\rho}{1-|y|^2}dy+\frac{C}{s^{2b}}\int_{B}w^2 \rho dy+ \frac{C}{s^{2b}}\int_{B} |\nabla w|^2 \rho (1-|y|^2)dy.\\
  \end{equation}
Combining $\eqref{112}$, $\eqref{6}$, $\eqref{10}$, $\eqref{12}$ and $\eqref{9},$ we write

\begin{eqnarray}
\frac{d}{ds}(J(w(s),s))&\leq &\frac{p+3}{2s^{b}}E(s) -\frac{p-1}{4s^{b}}\int_{B}(|\nabla w|^2-(y.\nabla w)^2)\rho dy
  -\frac{p+1}{2(p-1)s^{b}}\int_{B}w^2 \rho dy\nonumber\\
  &&-\frac{p-1}{2(p+1)s^{b}}\int_{B}|w|^{p+1}\rho dy+\frac{C}{s^{2b}}\int_{B}w^2 \rho dy
  + \frac{C}{s^{2b}}\int_{B}|\nabla w|^2(1-|y|^2)\rho dy\nonumber\\
  &&+\alpha \int_{B}(\partial_{s}w)^2\frac{\rho}{1-|y|^2}dy
+ C e^{\frac{-(p+1)s}{p-1}} +\frac{C}{s^{a+b}}\int_{B}|w|^{p+1}\rho dy.
\end{eqnarray}
Since $|y.\nabla w|\leq |y||\nabla w|$, it follows that
  $$\int_{B}|\nabla w|^2(1-|y|^2)\rho dy\leq \int_{B}(|\nabla w|^2-(y.\nabla w)^2)\rho dy.$$
  This leads finally to
  \begin{eqnarray*}
  \frac{d}{ds}(J(w(s),s))&\leq &\alpha \int_{B}(\partial_{s}w)^2\frac{\rho}{1-|y|^2}dy +\frac{p+3}{2s^{b}}E(s)-\frac{p-1}{4s^{b}}
\int_{B}|\nabla w|^2(1-|y|^2)\rho dy\\
&&-\frac{p+1}{2(p-1)s^{b}}\int_{B} w^2 \rho dy-\frac{p-1}{2(p+1)s^{b}}\int_{B}|w|^{p+1}\rho dy+\Sigma_{1}(s),
\end{eqnarray*}
where $\Sigma_{1}(s)$ satisfies the following inequality
$$\Sigma_{1}(s)\leq C e^{\frac{-(p+1)s}{p-1}}+  \frac{C}{s^{2b}}\int_{B} w^2 \rho dy+ \frac{C}{s^{2b}}\int_{B}|\nabla w|^2(1-|y|^2)\rho dy + \frac{C}{s^{a+b}}\int_{B}|w|^{p+1}\rho dy,$$
which ends the proof of Lemma 2.2.

\Box
\bigskip

For the reader's convenience we give the details of the proof of Theorem 1 in the following subsection.
\subsection{Proof of Theorem 1}
{\it Proof}: Before going into the proof, let's recall  from $\eqref{4}$ that
 $$E(w(s),s)=E_{0}(w(s))+I(w(s),s)+J(w(s),s).$$
 Now, according to Lemma 2.1 and Lemma 2.2, we have
 \begin{eqnarray*}
\frac{d}{ds}(E(w(s),s))&\leq &Ce^{\frac{-(p+1)s}{p-1}}+\frac{p+3}{2s^{b}}E(w(s),s)-\alpha \int_{B} (\partial_{s}w)^2 \frac{\rho}{1-|y|^2}dy\\
&&+ \Big(\frac{C}{s^{b}}-\frac{p-1}{4}\Big)\frac{1}{s^{b}}\int_{B} |\nabla w|^2(1-|y|^2)\rho dy\\
 &&+\Big(\frac{C}{s^{b}}-\frac{p+1}{2(p-1)}\Big)\frac{1}{s^{b}}\int_{B} w^2 \rho dy\\
 &&+\Big(\frac{C}{s^{a-b}}+\frac{C}{s^{a}}-\frac{(p-1)}{2(p+1)}\Big)\frac{1}{s^{b}}\int_{B}|w|^{p+1}\rho dy.
 \end{eqnarray*}
 Then, we consider $S_{1}\geq 1$ such that, for all $s\geq \max (S_{1},s_{0})$, we have:
 $$\frac{C}{s^{b}}-\frac{p-1}{4}\leq 0,\,\,\,\frac{C}{s^{b}}-\frac{p+1}{2(p-1)}\leq 0,\,\,\,\frac{C}{s^{a-b}}+\frac{C}{s^{a}}-\frac{(p-1)}{2(p+1)}\leq 0.$$
  Thus, implies that for all $s \geq \max (S_{1},s_{0})$,
  \begin{equation}\label{os}
  \frac{d}{ds}(E(w(s),s))\leq Ce^{\frac{-(p+1)s}{p-1}}+\frac{p+3}{2s^{b}} E(w(s),s)-\alpha \int_{B} (\partial_{s}w)^2 \frac{\rho}{1-|y|^2}dy.\\
\end{equation}
Recalling that,
$$H(w(s),s)=\exp\Big(\frac{p+3}{(a-1)s^{b-1}}\Big) E(w(s),s)+\theta e^{\frac{-(p+1)s}{p-1}},$$  
   we get from straightforward computations 
  \begin{eqnarray}\label{as}
  \frac{d}{ds}(H(w(s),s)) &=&-\frac{p+3}{2s^{b}}\exp\Big(\frac{p+3}{(a-1)s^{b-1}}\Big) E(w(s),s)\\
  &&+\exp\Big(\frac{p+3}{(a-1)s^{b-1}}\Big)\frac{d}{ds}(E(w(s),s))-\theta \frac{(p+1)}{p-1}e^{\frac{-(p+1)s}{p-1}}.\qquad\nonumber
\end{eqnarray}
Therefore, estimates $\eqref{os}$ and $\eqref{as}$ lead to the following crucial estimate:
 \begin{eqnarray*}
 \frac{d}{ds}(H(w(s),s))&\leq &\Big(C\exp\Big(\frac{p+3}{(a-1)s^{b-1}}\Big)-\theta \frac{(p+1)}{p-1}\Big)e^{\frac{-(p+1)s}{p-1}}\\
 &&-\alpha \exp\Big(\frac{p+3}{(a-1)s^{b-1}}\Big) \int_{B} (\partial_{s}w)^2 \frac{\rho}{1-|y|^2}dy.
\end{eqnarray*}
Since, we have $1\leq \exp\Big(\frac{p+3}{(a-1)s^{b-1}}\Big)\leq \exp\Big(\frac{p+3}{(a-1)}\Big)$, we deduce for all $s\geq \max (S_{1},s_{0})$,
$$\frac{d}{ds}(H(w(s),s))\leq \Big(C-\theta \frac{(p+1)}{p-1}\Big)e^{\frac{-(p+1)s}{p-1}}-\alpha \int_{B} (\partial_{s}w)^2 \frac{\rho}{1-|y|^2}dy.$$
We then choose $\theta$ large enough, so that $C-\theta \frac{(p+1)}{p-1}\leq 0$, which yields
$$\frac{d}{ds}(H(w(s),s))\leq -\alpha \int_{B} (\partial_{s}w)^2 \frac{\rho}{1-|y|^2}dy.$$
A simple integration between $s_{1}$ and $s_{2}$ ensures the result. This ends the proof of  the Theorem 1.

\Box
\bigskip

\begin{section}{Proof of Theorem 2 }
\end{section}
Throughout this section, we give a blow-up criterion in the w(y,s) variable and conclude the proof
of Theorem 2.
\subsection{A blow-up criterion in the w(y,s) variable}
We now claim the following lemma:
\begin{lem}
There exists $S_{2}\geq S_{1}$, such that for all $s_{0} \in {\R}$ and $w$ solution of equation
$\eqref{1}$ defined to the right of $s_{0}$, such that $\|w\|_{L^{p+1}(B)}$ is locally bounded, if $H(w(s_{3}),s_3)<0$ for some $s_{3} \geq \max(S_{2},s_{0})$, then $w$ cannot be defined for all $(y,s)\in B\times [s_{3}+1,+\infty)$.
 \end{lem}
\begin{rem}
{\rm Before going into the proof of Lemma 3.1, let's remark that if $w=w_{x_{0},T_{0}}$ defined from a solution of $\eqref{y}$ by $\eqref{2}$ and $x_{0}$ is
a non-characteristic point, then $\|w\|_{H^{1}(B)}$ is locally bounded and
so is $\|w\|_{L^{p+1}(B)}$ by sobolev's embedding.}
\end{rem}
{\it Proof}: The argument is the same as in the corresponding part in \cite{cf}. We sketch the proof for the reader's  convenience. Arguing by contradiction, we assume that there exists a solution $w$ on
$B$, defined for all $(y,s)\in B\times [s_{3}+1,+\infty)$, with $H(w(s_{3}),s_{3})<0$. Since the energy
$H(w(s_{3}),s_{3}))$ decreases in time, we have $H(w(s_{3}+1),s_{3}+1)<0$.
\\Consider now for $\delta>0$ the function $w^{\delta}(y,s)$ for $(y,s)\in B \times [s_{3}+1,+\infty)$,
defined for all $s\geq s_{3}+1$  and $y\in B$, by
\begin{equation}
w^{\delta}(y,s)=\frac{1}{(1+\delta e^{s})^{\frac{2}{p-1}}}w\Big(\frac{y}{1+\delta e^{s}},-\log(\delta+e^{-s})\Big),\,\,\,\,\,\,\,\,\\
\end{equation}
we have three observations:
\begin{itemize}
\item (A) Note that  $w^{\delta}$ is defined for all $(y,s)\in B\times [s_{3}+1,+\infty)$, whenever $\delta>0$ is small enough so that $-\log(\delta+e^{-s_{3}-1})\geq s_{3}$.
\item  (B) By construction, $w^{\delta}$ is also a solution of $\eqref{1}$ (indeed,
let $u$ be such that $w=w_{0,0}$ in definition $\eqref{2}$. Then $u$ is a solution of $\eqref{y}$ and $w^{\delta}=w_{-\delta,0}$
is defined as in $\eqref{2}$; so $w^{\delta}$ is also a solution of $\eqref{1}$).
\item (C) For $\delta$ small enough, we have $H(w^{\delta}(s_{3}+1),s_{3}+1)<0$ by continuity of the function $\delta\mapsto H(w^{\delta}(s_{3}+1,s_{3}+1))$. \\
\end{itemize}
Now, we fix $\delta=\delta_{0}>0$ such that (A), (B) and (C) hold. Let us note that we have
\begin{equation}\label{1919}
-\frac{1}{s^{b}}\int_{B}w^{\delta_{0}}\partial_{s}w^{\delta_{0}}\rho dy \geq -\frac{1}{4}\int_{B}(\partial_{s}w^{\delta_{0}})^2\rho dy-\frac{1}{s^{2b}}\int_{B}
(w^{\delta_{0}})^2\rho dy.
\end{equation}
According to the inequality  $\eqref{150}$, we obtain
\begin{equation}\label{est}
-e^{\frac{-2(p+1)s}{p-1}} \int_{B} F(e^{\frac{2s}{p-1}}w^{\delta_{0}})\rho dy \geq -C e^{\frac{-(p+1)s}{p-1}}-\frac{C}{s^{a}}\int_{B}|w^{\delta_{0}}|^{p+1}\rho dy .\\
\end{equation}
We recall that
$$E(w^{\delta_{0}}(s),s)=E_{0}(w^{\delta_{0}}(s))-e^{\frac{-2(p+1)s}{p-1}} \int_{B} F(e^{\frac{2s}{p-1}}w^{\delta_{0}})\rho dy-\frac{1}{s^{b}}\int_{B}w^{\delta_{0}}\partial_{s}w^{\delta_{0}}\rho dy.\,\,\,\,\,\,\,\,\,\,\,\,\,\,\,\,\,\,\,\,\,\,\,\,\,\,\,\,\,\,\,\,\,\,\,\,\,\,\,\,\,\,\,\,\,\,\,\,\,\,\,\,\,$$
Plugging the estimates $\eqref{1919}$ and $\eqref{est}$ together, we obtain
$$E(w^{\delta_{0}}(s),s) \geq  E_{0}(w^{\delta_{0}}(s))-C e^{\frac{-(p+1)s}{p-1}}-\frac{C}{s^{a}}\int_{B}|w^{\delta_{0}}|^{p+1}\rho dy-\frac{1}{4}\int_{B}(\partial_{s}w^{\delta_{0}})^2\rho dy-\frac{1}{s^{2b}}\int_{B}(w^{\delta_{0}})^2\rho dy.$$
By using the fact that
$$E_{0}(w^{\delta_{0}}(s))\geq \int_{B}\Big(\frac{1}{2}(\partial_{s}w^{\delta_{0}})^2+\frac{p+1}{(p-1)^2}(w^{\delta_{0}})^2-\frac{1}{p+1}|w^{\delta_{0}}|^{p+1}\Big)\rho dy,$$
it follows that
 \begin{eqnarray*}
E(w^{\delta_{0}}(s),s)&\geq &\int_{B}\Big(\frac{1}{2}(\partial_{s}w^{\delta_{0}})^2+\frac{p+1}{(p-1)^2}(w^{\delta_{0}})^2-\frac{1}{p+1}|w^{\delta_{0}}|^{p+1}\Big)\rho dy
 -C e^{\frac{-(p+1)s}{p-1}}\\
 &&-\frac{C}{s^{a}}\int_{B}|w^{\delta_{0}}|^{p+1}\rho dy-\frac{1}{4}\int_{B}(\partial_{s}w^{\delta_{0}})^2\rho dy-\frac{1}{s^{2b}}\int_{B}(w^{\delta_{0}})^2\rho dy.
\end{eqnarray*}
Hence, for any $s\geq s_{3}+1$,
  \begin{eqnarray*}
 E(w^{\delta_{0}}(s),s)&\geq &\frac{1}{4}\int_{B} (\partial_{s}w^{\delta_{0}})^2\rho dy +\Big(\frac{p+1}{(p-1)^2}-\frac{1}{s^{2b}}\Big)\int_{B}(w^{\delta_{0}})^2\rho dy\\
 &&-\Big( \frac{1}{p+1}+\frac{C}{s^{a}}\Big)\int_{B} |w^{\delta_{0}}|^{p+1}\rho dy-C e^{\frac{-(p+1)s}{p-1}}.
 \end{eqnarray*}
 We choose $s_{4}\geq s_{3}$ large enough, so that we have $\frac{p+1}{(p-1)^2}-\frac{1}{s^{2b}}\geq 0$. Then, we deduce, for all $s\geq s_{4}$,
 $$E(w^{\delta_{0}})\geq -\Big( \frac{1}{p+1}+\frac{C}{s^{a}}\Big)\int_{B} |w^{\delta_{0}}|^{p+1}\rho dy-C e^{\frac{-(p+1)s}{p-1}}.$$
By using the construction of $w^{\delta}$, we write
   $$E(w^{\delta_{0}}(s),s)\geq -\frac{\frac{1}{p+1}+\frac{C}{s^{a}}}{(1+\delta_{0} e^{s})^{\frac{2(p+1)}{p-1}}}\int_{B}|w\Big(\frac{y}{1+\delta_{0} e^{s}},-\log(\delta_{0}+ e^{-s})\Big)|^{p+1}\rho dy-C e^{\frac{-(p+1)s}{p-1}}.$$
  Since $\rho \leq 1$, the change of variable $z:=\frac{y}{1+\delta_{0} e^{s}},$ yields
$$ E(w^{\delta_{0}}(s),s)\geq  -\frac{ \frac{1}{p+1}+\frac{C}{s^{a}}}{(1+\delta_{0}e^{s})^{\frac{2(p+1)}{p-1}-N}}\int_{B}|w(z,-\log(\delta_{0}+e^{-s}))|^{p+1}dz-C e^{\frac{-(p+1)s}{p-1}}.$$
It is clear that $-\log(\delta_{0}+e^{-s})\rightarrow -\log(\delta_{0})$ as $s\rightarrow +\infty$ and since $\|w\|_{L^{p+1}(B)}$ is locally bounded
by hypothesis, by a continuity argument, it follows that the former integral remains bounded and
$$E(w^{\delta_{0}}(s),s)\geq-\frac{C}{(1+\delta_{0}e^{s})^{\frac{4}{p-1}+2-N}}-C e^{\frac{-(p+1)s}{p-1}} \longrightarrow 0,$$
as $s\longrightarrow +\infty$ (this is due to the fact that $\frac{4}{p-1}+2-N> 0$ which follows from the fact that $p<p_{c}$). So, thanks to $\eqref{3}$, it follows that
\begin{equation}\label{15}
\Liminf_{s \rightarrow +\infty} H(w^{\delta_{0}}(s),s)\geq 0.\\
\end{equation}
The inequality $\eqref{15}$ contradicts the inequality $H(w^{\delta_{0}}(s_{3}+1),s_{3}+1)<0$ (see item (C) above) and the fact that the energy $H$ decreases in time for $s\geq s_{3}$, which leads to the result.
This ends the proof of Lemma 3.1.
 \Box

 \subsection{Boundedness of the solution in similarity variables}

 We prove Theorem 2 here. Note that the lower bound follows from the finite speed
of propagation and wellposedness in $H^{1}\times L^{2}$. For a detailed argument in the similar case of equation $\eqref{16}$, see Lemma 3.1 (page 1136) in \cite{kl}.
\\We consider $u$ a solution of $\eqref{y}$ which is defined under the graph of $x\mapsto T(x)$
 and $x_{0}$ is a non-characteristic point. Given some $T_{0}\in (0,T(x_{0})]$, we introduce
$w_{x_{0},T_{0}}$ defined in $\eqref{2}$, and write $w$ for simplicity, when there is no ambiguity. We
aim at bounding $\|(w,\partial_{s}w)(s)\|_{H^{1}\times L^{2}}$ for $s$ large.
\\As in \cite{fh3}, by combining Theorem 1 and Lemma 3.1 (use in particular the remark after
that lemma) we get the following bounds:
\begin{cor}
{\rm For all} $s \geq \hat{s_{3}}=\hat{S_{3}}(T_{0})=\max(S_{3},-\log(T_{0}))$ , $s_{2}\geq s_{1} \geq  \hat{s_{3}}$,
{\rm it holds that}
$$-C \leq E(w(s),s)\leq M_{0}.$$
$$ \int_{s_{1}}^{s_{2}}\int_{B} (\partial_{s}w)^2 \frac{\rho}{1-|y|^2}dyds \leq M_{0}.$$
\end{cor}
Starting from these bounds, the proof of Theorem 2 is similar to the proof in \cite{fh3}, \cite{kl} and \cite{fh4}. To be more complete and in order to state our main result in a clear way, let us mention that the unique difference lays in the logarithmic term. In our opinion, handling these terms is straighforward
  in all the steps of the proof, except for the first step, where we bound the time averages of the $L^{p+1}_{\rho}(B)$ norm of $w$.
  For that reason, we only give that step and refer to \cite{fh3}, \cite{kl} and \cite{fh4}, for the remaining steps in the proof of Theorem 2.
   This is the step we prove here (in the following, $K_{3}$ denotes a constant that depends on $p$, $N$, $C$).
\bigskip
\begin{prop}
{\rm For all} $s\geq 1+\hat{s_{3}}$;
$$\int_{s}^{s+1}\int_{B} |w|^{p+1}\rho dy \leq K_{3}.$$
\end{prop}
{\it Proof}:
The proof of Proposition 3.3 is the same as in Hamza and Zaag \cite{MH} and \cite{MH1}. Exceptionally, the unique difference lays in the logarithmic term where we use the same technique as in the proof of Lemma 2.1 in Section 1.
Since the derivation of Theorem 2 from Proposition 3.3 is the same as in the non-perturbed case
treated in \cite{fh3}, \cite{kl} and \cite{fh4}, (up to some very minor changes), this concludes the proof of Theorem 2.
  \Box
 
 \no Acknowledgements: The authors are grateful
to Hatem Zaag for many fruitful discussions.



\noindent{\bf Address}:\\
Universit\'e de Tunis El Manar, Facult\'e des Sciences de Tunis, LR03ES04 \'Equations aux d\'eriv\'ees partielles et applications, 2092 Tunis, Tunisie\\
\vspace{-7mm}
\begin{verbatim}
e-mail: ma.hamza@fst.rnu.tn
e-mail: saidi.omar@hotmail.fr
\end{verbatim}

\end{document}